\documentclass[french,10pt,twoside,a4paper]{article}

\usepackage[latin1]{inputenc}
\usepackage[T1]{fontenc}
\usepackage{babel,indentfirst}
\usepackage{amsmath}
\usepackage{amssymb}

\newlength{\saut}
\newenvironment{disarray}
 {\everymath{\displaystyle\everymath{}}\array}
 {\endarray}

{\newenvironment{demo}{\begin{list}{}{%
}\item \noindent{\sc  D\'emonstration~-}~
}
{\vskip-\baselineskip\hfill $\blacksquare$\end{list}

\vspace{\saut}

}

\newtheorem{prop}{Proposition}[section]
\newtheorem{lem}[prop]{Lemme}
\newtheorem{cor}[prop]{Corollaire}

\newtheorem{thm}{Th\'eor\`eme}
\newtheorem{conj}{Conjecture}
\newtheorem{ppts}[prop]{Propriétés}
\newtheorem{defi}{Définition}[section]
\newtheorem{fact}[prop]{Fait}

\setlength{\marginparwidth}{30pt}
\setlength{\textwidth}{420pt}
\setlength{\evensidemargin}{20pt}
\setlength{\oddsidemargin}{20pt}

\DeclareMathOperator{\codim}{codim}
\DeclareMathOperator{\card}{Card}

\newcommand{\bg}{\mathbf b}
\DeclareMathOperator{\Ecc}{Ecc}
\newcommand{\Ac}{\mathcal A}
\newcommand{\Bc}{\mathcal B}

\newcommand{\Mc}{\mathcal M}
\newcommand{\Oc}{\mathcal O}
\newcommand{\Pc}{\mathcal P}

\newcommand{\ug}{\mathbf u}
\newcommand{\x}{\mathbf x}
\newcommand{\y}{\mathbf y}

\newcommand{\alphab}{\boldsymbol{\alpha}}
\newcommand{\betab}{\boldsymbol{\beta}}

\newcommand{\lambdab}{\boldsymbol{\lambda}}
\newcommand{\mub}{\boldsymbol{\mu}}

\newcommand{\thetab}{\boldsymbol{\theta}}

\newcommand{\bd}{\begin{demo}}
\newcommand{\ed}{\end{demo}}
\newcommand{\Cr}{{\mathbb C}}
\newcommand{\E}{{\mathbf E}}

\newcommand{\Fr}{{\mathbb F}}
\newcommand{\Lr}{{\mathbb L}}
\newcommand{\Nr}{{\mathbb N}}
\newcommand{\Prn}{{\mathbb P}}
\newcommand{\Qr}{{\mathbb Q}}

\newcommand{\Zr}{{\mathbb Z}}
\newcommand{\Qb}{\overline{\mathbb Q}}
\newcommand{\Qbnt}{\Qb^*\setminus\Qb_{\text tors}^*}
\newcommand{\Gm}{{\mathbb G}_m}
\newcommand{\io}{\omega}
\newcommand{\me}{ \hat{\mu}_{\text{ess}} }
\newcommand{\pg}{\bullet}

\title{
Minoration effective de la hauteur des points d'une courbe de $\Gm^2$ définie sur $\Qr$.
}
\date{}
\author{Corentin~{\sc Pontreau}.}

\begin{document}

\newpage

\setcounter{page}{0} 
\begin{center}
{\LARGE Minoration effective de la hauteur des points

\smallskip

d'une courbe de $\Gm^2$ définie sur $\Qr$.
}

\bigskip

{\large Corentin~{\sc Pontreau}.}
\end{center}

\smallskip

\thispagestyle{empty}

\begin{center}
\end{center}

\bigskip

\begin{abstract}

We are concerned here with Lehmer's problem in dimension $2$~; we give a lower bound for the height of a non-torsion point of $\Gm^2$ on a non-torsion curve defined over $\Qr$, depending on the degree of the curve only. We have first been inspired by~\cite{Am-Da3}; we develop a new approach, inherent in the dimension two (or more precisely the codimension two), and then obtain a better result where the error's term is improved significantly, moreover we give an explicit expression for the constant.
\end{abstract}

\bigskip

\maketitle


\section{Introduction.}

Dans tout cet article nous considérerons le plongement naturel de $\Gm^n:=\Gm^n(\Qb)$ dans $\Prn^n$, défini par $(\alpha_1,\ldots,\alpha_n)\mapsto(1:\alpha_1:\ldots:\alpha_n)$, en particulier la hauteur d'un point $\alphab$ sera la hauteur de Weil logarithmique du point projectif correspondant, soit, si $k$ est un corps de nombres contenant $\alpha_1,\dots,\alpha_n$~:
$$
h(\alphab):=\sum_{v\in\Mc_k} \frac{[k_v:\Qr_v]}{[k:\Qr]}\log\left(\max\{|\alpha_0|_v,\dots,|\alpha_n|_v\}\right).
$$
De même si $V$ est une sous-variété de $\Gm^n$, nous noterons $\deg(V)$ le degré de son adhérence de Zariski dans $\Prn^n$. Enfin, si $F$ est un polynôme à coefficients algébriques, nous noterons $h(F)$ la hauteur du point projectif défini par ses coefficients. 

Nous utiliserons la structure naturelle de groupe commutatif (donc de $\Zr$-module) de $\Gm^n$ ainsi, si $\alphab$ et $l$ désignent respectivement des éléments de $\Gm^n$ et de $\Zr$, alors, pour toute variété $V$, nous noterons $\alphab\cdot V:=\{\alphab\betab\ |\ \betab\in V\}$ et $[l]V:=\{\betab^l\ |\ \betab\in V\}$.

Rappelons que l'indice d'obstruction $\io_{\Qr}(\alphab)$ d'un point $\alphab\in\Gm^n$ n'est autre que le plus petit degré d'une hypersurface de $\Gm^n$ définie sur $\Qr$ passant par $\alphab$. Remarquons de plus que cette quantité est contr\^olée par le degré du corps de définition de $\alphab$~; en effet un argument 
d'algèbre linéaire nous donne l'inégalité~:
\begin{equation}
\label{deg_obs}
\io_{\Qr}(\alphab)\leq n(\deg_{\Qr}(\alphab))^{1/n}
\end{equation}

\noindent (qui est d'ailleurs une égalité si $n=1$). Dans~\cite{Am-Da} les auteurs proposent une conjecture  ({\it c.f.\ } Conjecture 1.3\footnote{Dans {\it op. cit.} l'indice d'obstruction $\io_{\Qr}(\alphab)$ est noté $\delta(\alphab)$.}) généralisant celle de {\sc D.H. Lehmer} ({\it c.f.} \cite{Le}) en dimension supérieure :

\begin{conj}\label{AD1conj}
Pour tout entier $n\geq 1$, il existe une constante $c(n)>0$ telle que, pour tout $\alphab\in\Gm^n$ dont les coordonnées sont multiplicativement indépendantes, on ait :
$$
h(\alphab)\geq \frac{c(n)}{\io_{\Qr}(\alphab)}.
$$
\end{conj}

Dans le même article, ils montrent que cette conjecture est vraie à un facteur <<$\log$>> près, généralisant un  théorème de {\sc E.~Dobrowolski} ({\it c.f.} \cite{Do}) en dimension supérieure. Ils poursuivent cette étude dans~\cite{Am-Da3} et affinent ce résultat~:

\newpage

\begin{thm}\label{AD3}
Soit $V$ une sous-variété algébrique de $\Gm^n$, définie sur une extension cyclotomique $K$ de $\Qr$, intersection d'hypersurfaces de $\Gm^n$ définies sur $K$ et de degré au plus $\io$.

Alors il existe une constante $c'(n)>0$ (effectivement calculable) telle que l'on ait, pour tout $\alphab\in V$ n'appartenant à aucune sous-variété de torsion\footnote{on appellera sous-variété de torsion une réunion de translatés de sous-groupes algébriques de $\Gm^n$ par des points de torsion.} de $V$~:
$$
h(\alphab)\geq \frac{c'(n)}{\io}\left(\log(3[K:\Qr]\io)\right)^{-\kappa(n)},
$$
où $\kappa(n):=2n(n+1)!^n-1$.
\end{thm}

Remarquons qu'un point $\alphab$ est à coordonnées multiplicativement indépendantes si et seulement s'il n'appartient à aucune sous-variété de torsion de $\Gm^n$, auquel cas le théorème~\ref{AD3} nous donne la minoration~:
$$
h(\alphab)\geq \frac{c(n)}{\io_{\Qr}(\alphab)}\left(\log(3\io_{\Qr}(\alphab))\right)^{-\kappa(n)}\ ;
$$
c'est précisément le résultat obtenu dans~\cite{Am-Da2}.

\bigskip

Notre principal résultat est le suivant, analogue du théorème~\ref{AD3} où l'exposant du $\log$ ($\kappa(2)=143$) est sensiblement amélioré. De plus, contrairement à ce dernier, il est complètement explicite~:

\noindent\rule{\textwidth}{.2mm}

\vspace{-2mm}

\begin{thm}\label{Mthm}
Soit $V$ une courbe de $\Gm^2$ définie sur $\Qr$ et $\Qr$-irréductible de degré $\io$ qui ne soit pas de torsion, et soit $\alphab\in V\setminus(\Gm^2)_{\text{tors}}$, on a~:
$$
h(\alphab) \geq \frac{1,2\cdot 10^{-16}}{\io} \frac{\big(\log\log\io'\big)^{11}}{\big(\log\io'\big)^{13}}
$$
où $\io':=\max\{\io,16\}$.
\end{thm}

\vspace{-2mm}

\noindent\rule{\textwidth}{.2mm}

Notons qu'ici l'hypothèse <<$\alphab\in V$ n'appartenant à aucune sous-variété de torsion de $V$>> se réduit à <<$\alphab$ non de torsion>>, hypothèse minimale puisque les points de torsion sont précisément les points de hauteur nulle.

\bigskip

F. {\sc Amoroso} et S. {\sc David} montrent dans~\cite{Am-Da} que si $\alphab\in\Gm^n$ est un point à coordonnées multiplicativement indépendantes, alors il existe une constante $C(n)$ telle que :
\begin{equation}\label{AmDaCor}
\left(h(\alpha_1)\dots h(\alpha_n)\right)^{1/n}\geq \frac{C(n)^{1/n}}{D^{1/n}} (\log 3D)^{-\kappa(n)},
\end{equation}
où $D:=[\Qr(\alphab):\Qr]$ et $\kappa(n)$ est la même quantité que dans le théorème~\ref{AD3}. Une observation de Bilu (voir~\cite{Bi}) mène à la minoration :
$$
\left(h(\alpha_1)\dots h(\alpha_n)\right)^{1/n}\geq \frac{C(2)^{1/2}}{D^{1/2}} (\log 3D)^{-\kappa(2)}.
$$
En effet, si l'on applique l'inégalité~(\ref{AmDaCor}) en dimension 2 pour $(\alpha_1,\alpha_2),\ldots,(\alpha_{n-1},\alpha_n)$ et $(\alpha_n,\alpha_1)$, on obtient :
$$
\left(\prod_{i=1}^n h(\alpha_i)\right)^{1/n}\geq \frac{C(2)^{1/2}}{D^{1/2}} (\log 3D)^{-\kappa(2)}.
$$
Le théorème~\ref{Mthm} nous permet d'expliciter la constante $C(2)$ et d'améliorer l'exposant $\kappa(2)$~:

\begin{flushleft}
\begin{minipage}{\textwidth}
\noindent\rule{\textwidth}{.2mm}

\vspace{-1mm}

\begin{cor}\label{Mcor}
Soient $\alphab\in\Gm^n$ et $\sigma$ une permutation de $\{1,\dots,n\}$ telle que, pour tout $i$ dans $\{1,\ldots,n\}$, les coordonnées $\alpha_i$ et $\alpha_{\sigma(i)}$ soient multiplicativement indépendantes, alors :
$$
\left(h(\alpha_1)\dots h(\alpha_n)\right)^{1/n}\geq \frac{2\cdot 10^{-21}}{D^{1/2}} (\log 3D)^{-13},
$$
où $D:=[\Qr(\alphab):\Qr]$.
\end{cor}

\vspace{-3mm}

\noindent\rule{\textwidth}{.2mm}
\end{minipage}

\end{flushleft}

Dès que l'on fixe un entier $n\geq 3$, ce résultat est moins bon que~(\ref{AmDaCor}) lorsque $D$ est grand, néanmoins il est bien meilleur pour les petites valeurs de $D$ (rappelons que $\kappa(n):=2n(n+1)!^n-1$), de plus il est entièrement explicite et les conditions sur $\alphab$ sont un peu plus faibles.

\bigskip

\noindent{\bf Remerciements.} Je tiens à remercier Francesco Amoroso pour l'aide et le soutien qu'il m'a apportés tout au long de ce travail. Je voudrais également tout particulièrement remercier Federico Pellarin pour les lectures très soignées qu'il a pu faire sur des versions préliminaires de ce travail, ainsi que les nombreux conseils et remarques qu'il a pu me donner par la suite à ce sujet.

\section{Schéma de la preuve}

Dans le paragraphe~\ref{parResPrel} on développe les outils d'une démonstration de transcendance (lemme de Siegel, extrapolation) et au paragraphe~\ref{parVE}, on montre des minorations explicites pour les courbes non de torsion de $\Gm^2$ et des points de $\Gm$. La démonstration du théorème~\ref{Mthm} à proprement parler est l'objet du paragraphe~\ref{DemoMthm}. La stratégie est la suivante~: par l'absurde on suppose la hauteur de $\alphab$ petite, on peut alors construire un polynôme s'annulant sur $V$ avec multiplicité, de degré et de hauteur contrôlés via un lemme de Siegel (proposition~\ref{Siegel2}). Ensuite on extrapole dans le paragraphe~\ref{MthmExrtap} en montrant, grâce à la proposition~\ref{premiers}, que ce polynôme s'annule sur les puissances $\alphab^{pq}$, où $p$ et $q$ parcourent des ensembles de premiers $\Pc_1$ et $\Pc_2$.

Nous reprenons au paragraphe~\ref{sub5.3} le lemme de zéros de~\cite{Am-Da3} (théorème~2.6)~; nous obtenons ainsi une suite décroissante $Y_1\supseteq Y_2\supseteq Y_3$ de sous-variétés de $\Gm^2$ contenant des puissances $\alphab^{pq}$ de $\alphab$. Deux de ces variétés étant de même dimension, on obtient une sous-variété obstructrice $Z$, composante $\Qr$-irréductible de $Y_3$ ou $Y_2$ contenant une puissance $\alphab^\ell$ de $\alphab$ et dont on contrôle le degré. Deux cas se présentent alors.

Si la variété obstructrice $Z$ est de dimension $0$ (paragraphe~\ref{CasDim0}), alors $Z$ est simplement l'union des conjugués de $\alphab$. Comme ici $\ell=1$, on arrive, grâce notamment à l'inégalité~(\ref{deg_obs}), à un encadrement du type~:
$$
\card(\Pc_2)\io_{\Qr}(\alphab)^2\ll\card(\Pc_2)\deg(Z)\ll (\log\io_{\Qr}(\alphab))^a\io_{\Qr}(\alphab)^2
$$
ainsi, de par nos choix de paramètres, on arrive à une contradiction. 

Dans le cas où la variété obstructrice $Z$ est de dimension $1$ (paragraphe~\ref{CasDim1}), la puissance $\ell$ est a priori différente de $1$. On peut obtenir l'encadrement~:
$$
\card(\Pc_1)\io_{\Qr}(\alphab^\ell)\ll\card(\Pc_1)\deg(Z)\ll (\log\io_{\Qr}(\alphab))^b\io_{\Qr}(\alphab)
$$
ainsi l'indice d'obstruction de $\alphab^\ell$ est très petit par rapport à celui $\alphab$. Ceci n'étant pas suffisant pour conclure, dans~\cite{Am-Da3}  les auteurs utilisent un argument de descente pour arriver à une contradiction. Ici notre démonstration diffère~; des arguments plus simples nous donnent de meilleurs résultats. On travaille avec la hauteur normalisée~; on majore celle de $Z$ en fonction de la hauteur de notre fonction auxiliaire $F$, sur laquelle on a un bon contrôle~:
$$
\card(\Pc_1)^2\hat{h}(Z)\ll h(F).
$$
Si $Z$ n'est pas une courbe de torsion, on arrive à une contradiction en utilisant une minoration explicite de $\hat{h}(Z)$ (proposition~\ref{dim2}). Dans le cas contraire, on se ramène dans le lemme~\ref{torsion} à une étude en dimension $1$, auquel cas, {\em via} la proposition~\ref{the:Dobrowolski2}, on obtient une minoration de $h(\alphab)$.

\bigskip

\section{Résultats préliminaires.}\label{parResPrel}

\subsection{Premiers exceptionnels.}

Dans ce paragraphe $V$ désigne une sous-variété de $\Gm^n$ irréductible sur son corps de définition. Nous allons voir que pour tout nombre premier $p\ $ sauf pour certains appartenant à un ensemble exceptionnel $\Ecc(V)$, introduit dans~\cite{Am-Da}, la variété $[p]V$ a un bon comportement, dans un sens que nous allons préciser.

\bigskip

\begin{defi}\label{defEcc}
On note $W_1,\ldots,W_k$ les composantes $\Qb$-irréductibles de $V$, on pose :
$$
\begin{array}{lcl}
\Ecc(V) & := & \big\{l\in\Zr\ |\ \exists i,j,\ i<j\ ;\  [l](W_i)=[l]W_j)\big\}\\[1mm]
        &     & \bigcup \big\{l\in\Zr\ |\  \exists i\ ;\ \deg([l]W_i)<\deg(W_i)\big\}.
\end{array}
$$ 
\end{defi}

\bigskip

La proposition 2.4 de~\cite{Am-Da} donne des informations sur cet ensemble~; nous en rappelons quelques propriétés.

\begin{prop}\label{CardEcc} Nous avons
$$
\card\left(\Ecc(V)\cap\{p\ \text{premier}\}\right)\leq \frac{\dim(V)+1}{\log 2} \log \big(\deg(V)\big).
$$
De plus, si $\Lambda$ est un ensemble fini de nombres premiers et si $V$ n'est pas de torsion, alors :
$$
\deg\Big(\bigcup_{p\in\Lambda} [p]V\Big) \geq \card\left(\Lambda\setminus\Ecc(V)\right)\times\deg(V).
$$
\end{prop}

\bigskip

Nous utiliserons dans la suite cette proposition dans le cas où $\Lambda$ est l'ensemble des premiers dans $[N/2,N]$, pour un certain paramètre $N$, d'où la nécessité du lemme suivant :

\begin{lem}\label{PiN}
Pour tout réel $x$ on note $\pi(x)$ le nombre de premiers inférieurs ou égaux à $x$. Pour $N\geq 41$ on a
$$
\pi(N)-\pi(N/2)\geq 0,41\frac{N}{\log N}.
$$
\end{lem}

\bd Le théorème~1 de \cite{Ro-Sc} nous donne :
$$
\forall x\geq 59,\quad \frac{x}{\log x}+\frac{3x}{2(\log x)^2}\ \geq\ \pi(x) > \frac{x}{\log x}+\frac{x}{2(\log x)^2}, 
$$
si on note $c_N:=\log (N/2)/\log (N)$ on en déduit :
$$
\begin{disarray}{lcl}\label{piN}
\pi(N)-\pi(N/2) & > & 
\frac{N}{\log N}+\frac{N}{2(\log N)^2}-\left(\frac{N}{2c_N\log N}+\frac{3N}{4(c_N\log N)^2}\right)\\[4mm]
                                                       & = &
\frac{N}{\log N}\left(1-\frac{1}{2c_N}-\Big(\frac{3}{4c_N^2}-\frac{1}{2}\Big)\frac{1}{\log N}\right)
\end{disarray}
$$
Ainsi, pour $N\geq 5000$, nous avons bien l'inégalité voulue et une vérification numérique pour les petites valeurs de $N$ nous permet de conclure.
\ed


\subsection{Construction de la fonction auxiliaire.}

On cherche ici un polynôme $F$ de degr\'e $\leq L$ nul en $\alphab$ \`a un ordre $\geq T$. Fixons dans un premier temps quelques notations, $k$ désignant un corps de nombres.

\begin{itemize}
\item On notera $k[\x]:=k[x_1,\dots,x_n]$.
\item Pour $\mub,\lambdab\in\Nr^n$, on pose :
$$
\begin{disarray}{lcl}
{\mub \choose \lambdab} & := & {\mu_1 \choose \lambda_1}\cdots{\mu_n \choose \lambda_n}.\\[3mm]
D_{\lambdab} & := & \frac{1}{\lambdab!}\frac{\partial^{\lambdab}}{\partial \x^{\lambdab}}=\frac{1}{\lambda_1!\cdots\lambda_n!}\frac{\partial^{\lambda_1}}{\partial x_1^{\lambda_1}}\cdots\frac{\partial^{\lambda_n}}{\partial x_2^{\lambda_n}}.
\end{disarray}
$$
\item On dira que $\alphab\in\Qb^n$ est racine de $F\in\Cr[\x]$ de multiplicité au moins $T\in\Nr$ si 
$$
\forall\lambdab\in\Nr^n,\ \lambda_1+\ldots+\lambda_n=|\lambdab|<T\ \Longrightarrow\ D_{\lambdab} (F)(\alphab)=0.
$$
\item On notera $\io_k(\alphab):=\min \big\{ \deg(F)\ |\ F\in k[\x]\setminus\{0\}, F(\alphab)=0 \big\}$ et $k[\x]_{\leq L}$ le $k$-espace vectoriel des polynômes de degré total $\leq L$.
\item Pour toute partie $S$ de $\Cr^n$ on pose : 
$$
E_k(S,L,T)  := \big\{ F\in k[\x]_{\leq L}\ |\ \forall\betab\in S,\ F \text{ nul en } \betab \text{ avec multiplicité } \geq T\big\}\footnotemark[1].
$$
\footnotetext[1]{
on a donc $E_k(S,L,T)=[\mathfrak{P}^{(T)}]_L\cup\{0\}$ si $S$ est une variété et $\mathfrak{P}$ est l'idéal de définition sur $k$ de sa clôture projective.
}
\end{itemize}

Nous aurons besoin dans la suite d'encadrements pour l'indice d'obstruction $\io_k(\alphab)$ et pour la dimension du $k$-espace vectoriel $E_k(\{\alphab\},L,T)$.

\bigskip

\begin{ppts}\label{pptdeg}
Soient $n\in\Nr^*$, $\alphab\in\Gm^n$ et $L,T\in\Nr$, on a :
\begin{enumerate}
\item $\dim E_k(\{\alphab\},L,T)\geq {L-T\io_k(\alphab)+ n \choose n}$.
\item $1\leq\io_\Qr(\alphab)\leq n[k(\alphab):k]^{1/n}$.
\end{enumerate}
\end{ppts}

\bd
\begin{enumerate}
\item Soit $F\in k[\x]$ non nul de degré $\omega_k(\alphab)$ tel que $F(\alphab)=0$. Pour tout $H\in k[\x]$ de degré inférieur ou égal à $L-T\omega_k(\alphab)$, on a $F^TH\in E_k(S,L,T)$, de plus le sous-espace vectoriel de $k[\x]$ des polynômes de degré $\leq L-T\omega_k(\alphab)$ est de dimension ${L-T\omega_k(\alphab)+n \choose n}$, d'où le résultat.

\item Posons $\omega:=[n[k(\alphab):k]^{1/n}]$  et considérons l'application linéaire :
$$
\begin{array}{ccc}
k[\x]_{\leq \omega} & \longrightarrow & k(\alphab)\\
        P                 & \longmapsto    & P(\alphab).
\end{array}
$$
Remarquons que $\dim_{ k} k[\x]_{\leq \omega}={\omega+n\choose n} \geq n^{-n}(\omega+1)^n$. Or $n^{-n}(\omega+1)^n>[k(\alphab):k]$, cette application n'est donc pas injective, {\it i.e.} il existe $P\in k[\x]_{\leq \omega}$ non nul tel que $P(\alphab)=0$, donc $\omega_k(\alphab)\leq \io\leq n[k(\alphab):k]^{1/n}$.
\end{enumerate}
\ed

Ci-dessous nous donnons la version du lemme de Siegel que nous utiliserons dans la suite (analogue à la proposition 2.1 de~\cite{Am-Da3}). 

\begin{prop}\label{Siegel2}
Soient $\theta$ un réel $>0$ et $\E\subset\{\alphab\in\Gm^n\ |\ h(\alphab)\leq\theta\}$ un ensemble fini non vide. Soient $k$ un corps de nombres et $L,T\in\Nr$. Si $E_k(\E,L,T)$ est non réduit à $\{0\}$, alors il existe un polynôme $F\in E_k(\E,L,T)\cap \Oc_k[\x]$ non nul tel que~:
\begin{equation}\label{IneqSiegel2}
h(F)\leq \frac{r}{N-r}\Big( (T+n-1)\log (L+1)+L\theta\Big)+\log c_k.
\end{equation}
où $c_k:=\left\{(\frac{2}{\pi})^s\sqrt{|\Delta_k|}\right\}^{\frac{1}{[k:\Qr]}}$, $s$ le nombre de places complexes de $k$ et $\Delta_k$ son discriminant, $N:=\dim_kk[\x]_{\leq L}$ et $r:=\dim_kk[\x]_{\leq L}-\dim_k E_k(\E,L,T)$.\\
\end{prop}

\bd
On reprend ici principalement les preuves de la proposition 4.2 de~\cite{Am-Da} et du théorème 7 de ~\cite{St-Va}. Fixons un ordre sur $\E$ et sur $\Nr^n$ et considérons la matrice $A$ de taille $\card(\E)\cdot{T+n-1 \choose n}\times {L+n \choose n}$
d\'efinie par
$$
A:=\left( {\mub \choose \lambdab} \alphab^{\mub-\lambdab} \right)
$$
o\`u les lignes (respectivement les colonnes) sont index\'ees par les couples
$(\alphab,\lambdab)$, o\`u $\alphab\in \E$ et  $\lambdab\in\Nr^n$ est tel que $|\lambdab|\leq T-1$ (respectivement par les multi-indices $\mub\in\Nr^n$ tels que $|\mub|\leq L$). Autrement dit, si on pose
$$
\Ac:=\{ \x\in k^N,\ A\x=0\}\ ,
$$
alors $\Ac=E_k(\E,L,T)$ d'où $r=\hbox{rang}(A)$. Soit $Y$ une matrice $N\times (N-r)$ à coefficients dans $k$ telle que
$\Ac$ soit l'image de l'application $k$-linéaire définie par $Y$. Comme $E_k(\E,L,T)$ est non réduit à $\{0\}$, on a rang$(Y)=N-r<N\ $; le théorème 8 de~\cite{Bo-Va} appliqué à $Y$ montre alors qu'il existe $N-r$ vecteurs linéairement indépendants
$\ug_1,\ldots,\ug_{N-r}$ de $k^{N-r}$ tels que, si l'on pose $F_i:=Y\ug_i$, pour
$i=1,\ldots, N-r$, on ait $F_i\in \Oc_k^N$ pour tout $i$ et :
$$
\sum_{j=1}^{N-r} h(F_j) \leq \log H(Y)+(N-r)\log c_k= \log H(\mathcal{A})+(N-r)\log c_k\ ,
$$
o\`u $H(\mathcal{A})$ est la hauteur (non logarithmique) du sous-espace $\mathcal{A}$ et $H(Y)$ la hauteur du sous-espace engendré par ses lignes (définies p.~499 de~\cite{St-Va}).

Remarquons que $(F_1,\ldots,F_{N-r})$ est une une base de $\mathcal{A}$, en particulier il existe $F$ dans $\Ac\cap\Oc_k^N$ non nul tel que
\begin{equation}\label{choixPoly}
(N-r)h(F) \leq \log H(\mathcal{A})+(N-r)\log c_k,
\end{equation}
nous allons montrer que ce $F$ vérifie bien~(\ref{IneqSiegel2}).

\medskip

Soit $\Fr$ la clôture galoisienne de $k(\E)/k$, considérons la matrice :
$$
B:=\left(
\begin{array}{c}
\sigma_1A\\
\vdots\\
\sigma_RA
\end{array}
\right)
$$
où les $\sigma_i$ sont les éléments de Gal$(\Fr/k)$, et posons :
$$
\Bc:=\{\y\in\Fr\ |\ B\y=0\}.
$$
On a alors $\dim_{\Fr}(\Bc)=\dim_k(\Ac)$ et $H(\Ac)=H(\Bc)$ (voir~\cite{St-Va}, (2.31) page 506).

Soit $\tilde{B}$ une sous-matrice de $B$ de rang maximal ($\tilde{B}$ est une matrice $r\times
{L+n\choose n}$ de rang $r$), par le principe de dualité, (voir \cite{St-Va} p. 500, (2.2)), on a :
$$
H(\Ac)=H(\Bc)=H(\tilde{B}).
$$
En majorant $H(\tilde{B})$ par le produit des hauteurs de ses lignes (in\'egalit\'e de Hadamard, voir
\cite{Bo-Va}, \'equation (2.6)), on obtient :
\begin{equation}\label{majHautB}
\log H(\tilde{B}) \leq r \log \max \Big\{H(\bg^{(\alphab,\lambdab)}) \ |\ \alphab\in\E\ \text{et}\ |\lambdab|\leq T-1 \Big\},
\end{equation}
o\`u les $\bg^{(\alphab,\lambdab)}$ désignent les lignes de $\tilde{B}$~:
$$
\bg^{(\alphab,\lambdab)}=(\bg_{\mub}^{(\alphab,\lambdab)})_{|\mub|\leq L}=\bigg( {\mub \choose \lambdab} \alphab^{\mub-\lambdab} \bigg)_{|\mub|\leq L}.
$$
Soit $(\alphab,\lambdab)$ un multi-indice r\'ealisant ce maximum ($|\lambdab|\leq T-1$), on a :
$$
\begin{disarray}{r c l}
\Big( \sum_{|\mub|\leq L} {\mub \choose \lambdab}^2\Big)^{\frac{1}{2}} & \leq &
 \sum_{|\mub|\leq L} {\mub \choose \lambdab}
= \sum_{\mu_1=1}^{L} \cdots\sum_{\mu_n=1}^{L} {\mu_1 \choose \lambda_1} {\mu_n \choose \lambda_n}\\[.6cm]
& = & {L+1 \choose \lambda_1+1} \cdots {L+1 \choose \lambda_n+1}\leq (L+1)^{T+n-1}
\end{disarray}
$$
où l'on a utilisé $\sum_{\mu=1}^{L} {\mu \choose \lambda}= {L+1 \choose
\lambda+1}\leq (L+1)^{\lambda+1}.$

Notons $d$ le degré de $\Fr$ sur $\Qr$, en utilisant cette in\'egalit\'e nous trouvons, pour toute place archim\'edienne
$v\in \mathcal{M}_k$, 
$$
H_v(\bg^{(\alphab,\lambdab)})^{d}=\Big( \sum_{|\mub|\leq L}
|\bg^{(\alphab,\lambdab)}_{\mub}|_v^2\Big)^{ \frac{d_v}{2} }
\leq
(L+1)^{(T+n-1)d_v}\max\big\{1,|\alpha_1|_v,\ldots,|\alpha_n|_v\big\}
^{Ld_v}.
$$
Pour $v\in \mathcal{M}_k$ ultramétrique, on obtient :
$$
H_v(\bg^{(\alphab,\lambdab)})^{d}=\max_{|\mub|\leq L}
|\bg^{(\alphab,\lambdab)}_{\mub}|_v^{d_v} \leq
\max\big\{1,|\alpha_1|_v,\ldots,|\alpha_n|_v\big\}^{Ld_v}.
$$
En faisant le produit sur toutes les places on obtient :
$$
\begin{disarray}{r @{\ \leq\ } l}
H(\bg^{(\alphab,\lambdab)})^{d} &
 (L+1)^{(T+n-1)\sum_{v|\infty}d_v} H(\alphab)^{L d}\\
d\log \big(H(\bg^{(\alphab,\lambdab)})\big) &
d(T+n-1)\log(L+1) +d L h(\alphab)\ ,
\end{disarray}
$$
d'où, en reprenant l'inégalité~(\ref{majHautB})
$$
\log H(\tilde{B})\leq r\Big((T+n-1)\log(L+1) +L\theta\Big)\ ,
$$
donc~(\ref{choixPoly}) devient :
$$
\begin{disarray}{lcl}
h(F)    & \leq &  \frac{1}{N-r}\log H(\tilde{B})+\log c_k\\
        & \leq & \frac{r}{N-r} \cdot \Big( (T+n-1)\log(L+1) + L\theta \Big)+\log c_k.
\end{disarray}
$$
\ed

\bigskip

Dans la suite, on utilisera cette proposition dans le cas $n=2$ :

\begin{cor}\label{f_auxCor}
Soient $\alphab\in\Gm^2$, $T\in\Nr^*$, $D$ le degré de $\Qr(\alphab)$ sur $\Qr$, et $\io\geq\io_{\Qr}(\alphab)$.\\
Si $L=\min\{\ 2\io T^2,\ \left[(TD)^{1/2}(T+1)\right]\}$, alors il existe un polynôme $F\in E_{\Qr}(\{\alphab\},L,T)\cap \Zr[\x]$ non nul tel que :
$$
h(F)\leq \frac{1}{T-1}\big( (T+1)\log (L+1)+Lh(\alphab)\big).
$$
\end{cor}

\bd
Comme $D^{1/2}\geq \frac{1}{2}\io_{\Qr}(\alphab)$, on a $L\geq \frac{1}{2}\io_{\Qr}(\alphab)(T+1)T^{1/2}\geq \io_{\Qr}(\alphab)T$, en particulier $E_{\Qr}(\{\alphab\},L,T)$ n'est pas réduit à $\{0\}$, d'après la proposition~\ref{pptdeg}. Considérons le polynôme $F$ donné par la proposition~\ref{Siegel2}, on a :
$$
h(F) \leq \frac{r}{N-r}\big( (T+1)\log (L+1)+Lh(\alphab)\big)\ ,
$$
où $r:=\dim\Qr[\x]_{\leq L}-\dim_{\Qr} E_{\Qr}(\{\alphab\},L,T)$ et $N:=\dim\Qr[\x]_{\leq L}$.

On a, si $L=2\io T^2$ :
$$
\begin{disarray}{lcl}
\frac{r}{N-r}   &  \leq & \frac{{L+2 \choose 2}}{{L-\io T+2 \choose 2}}-1 
= \frac{L+1}{L-\io T+1}\times\frac{L+2}{L-\io T+2}-1\\[5mm]
        & \leq & \left(\frac{L}{L-\io T}\right)^2-1 
= \left(\frac{2T}{2T-1}\right)^2-1 \leq \frac{1}{T-1}.
\end{disarray}
$$
Sinon, si $L=\left[(TD)^{1/2}(T+1)\right]$~:
$$
\begin{disarray}{lcl}
\frac{r}{N-r}   &  \leq & \frac{{T+1 \choose 2}D}{{L+2 \choose 2}-{T+1 \choose 2}D} 
= \frac{T(T+1)D}{(L+1)(L+2)-T(T+1)D}\\[5mm]
        & \leq & \frac{T(T+1)D}{(T+1)^2TD-T(T+1)D} 
= \frac{1}{T} \leq \frac{1}{T-1}.
\end{disarray}
$$
\ed

\subsection{Extrapolation}
 Nous allons utiliser le lemme clef de Dobrowolski ({\it c.f.}~\cite{Do}) dans le cadre plus large de polynômes à plusieurs variables à coefficients dans un anneau d'entiers d'une extension cyclotomique de $\Qr$.
Soient $k/\Qr$ une extension galoisienne, $p$ un nombre premier non ramifié dans $k$ et $Q$ un idéal premier de $\Oc_k$ tel que $Q|p$. Si l'extension $k/\Qr$ est abélienne, alors l'automorphisme de Frobénius associé $\phi_{Q,p}\in$ Gal$(k/\Qr)$ ne dépend que de $p$ et on le notera $\phi_p$~; on a
$$
\forall\alpha\in\Oc_k,\ \phi_p(\alpha)\equiv \alpha^p\mod p\Oc_k.
$$

Dans tout ce paragraphe, on supposera $k/\Qr$ {\bf cyclotomique} et on notera $\Delta_k$ son discriminant. De plus $\alphab$ désignera un élément de $\Gm^n$, $\Fr$ la clôture galoisienne de $k(\alphab)$, et $L,T$ deux entiers naturels. Le résultat qui suit correspond au théorème 2.2 de~\cite{Am-Da3}~:

\begin{thm}\label{clef} 
Soit $F\in E_k(\{\alphab\},L,T)\cap\Oc_k[\x]$~; pour tout nombre premier $p\nmid\Delta_k$
et pour tout $v\in \mathcal{M}_{\Fr}$ divisant $p$, on a la majoration
$$
|F^{\phi_p}(\alphab^p)|_{v}\leq p^{-T}|F|_{v} \max
\{1,|\alpha_1|_v,\ldots,|\alpha_n|_v\}^{pL},
$$
o\`u l'on a not\'e $\alphab^p=(\alpha_1^p,\ldots,\alpha_n^p).$
\end{thm}

\begin{prop}\label{premiers}
Soit $F\in E_k(\{\alphab\},L,T)\cap\Oc_k[\x]$. Pour tout nombre premier $p\nmid\Delta_k$, le polynôme $F^{\phi_p}$ est nul en $\alphab^p$ à un ordre $T_1$ vérifiant :
$$
T_1\Big(\log (L+1)+\log p\Big)\ >\ T\log p-h(F)-pLh(\alphab)-n\log(L+1).
$$
\end{prop}

\bd
Soit $\lambdab\in\Nr^n$ tel que $|\lambdab|=T_1$ et $D_{\lambdab}(F)^{\phi_p}(\alphab^p)\neq 0$ (on peut supposer $T_1< T$). Soit $v\in\Mc_\Fr$~; on déduit de l'inégalité 
$$
\sum_{|\mub|\leq L} {\mub \choose \lambdab}\leq {L+1 \choose \lambda_1+1}\cdots {L+1 \choose \lambda_n+1}\leq (L+1)^{|\lambdab|+n}
$$
et de l'inégalité ultramétrique les majorations~:
$$
|D_{\lambdab}(F)^{\phi_p}(\alphab^p)|_v \leq
\left\{
\begin{array}{l}
|F|_v \cdot \max\{1,|\alpha_1|_v,\ldots,|\alpha_n|_v\}^{pL} \hfill \text{ si } v\nmid\infty \text{ et } v\nmid p\\[2mm]
(L+1)^{|\lambdab|+n} |F|_v \cdot \max\{1,|\alpha_1|_v,\ldots,|\alpha_n|_v\}^{pL} \text{ si } v\mid\infty
\end{array}
\right.
$$
De plus, si $v\mid p$, le théorème~\ref{clef} donne :
$$
|D_{\lambdab}(F)^{\phi_p}(\alphab^p)|_v\leq p^{-(T-|\lambdab|)}|F|_v \max\{1,|\alpha_1|_v,\ldots,|\alpha_n|_v\}^{pL}.
$$
On a, par la formule du produit :
$$
1=\prod_{v\in\mathcal{M}_\Fr}
|D_{\lambdab}(F)^{\phi_p}(\alphab^p)|_v^{[\Fr_v:k_v]/[\Fr:k]}.
$$
En passant au $\log$, et en utilisant les trois majorations obtenues ci-dessus on obtient~:
$$
0\leq (|\lambdab|+n)\log(L+1)+h(F)+pLh(\alphab)-(T-|\lambdab|)\log p\ ,
$$
soit
$$
|\lambdab|\Big(\log (L+1)+\log p\Big)\ >\ T\log p-h(F)-pLh(\alphab)-n\log (L+1).
$$
\ed

\section{Versions explicites de certaines minorations}\label{parVE}

\subsection{Une minoration pour les courbes}

Dans~\cite{Am-Da2}, F. Amoroso et S. David obtiennent une minoration de la hauteur d'une hypersurface de $\Gm^n$ définie sur $\Qr$ et $\Qr$-irr\'eductible qui n'est pas de torsion ; de plus notre résultat principal dans~\cite{DEA} donne une version explicite de ce résultat. Nous reprenons celui-ci en améliorant la constante pour $n=2$, cas qui nous intéresse ici.

\begin{prop}\label{dim2}
Soit $V$ une courbe définie sur $\Qr$ et $\Qr$-irr\'eductible de $\Gm^2$ de degr\'e $D$. Alors, si $V$ n'est pas de torsion, on a
$$
\hat{h}(V)\geq 5^{-6} \left(\frac{\log\log D'}{\log D'}\right)^3\ ,
$$
où\footnote{Nous avons choisi de mettre $D'$ dans le $\log$ car la fonction $x\mapsto \frac{\log(x)}{\log\log(x)}$ est croissante sur $[16,+\infty[$, et afin de pouvoir minorer $\log\log D'$ par $1$ dans les calculs.} $D':=\max\{16,D\}$.
\end{prop}

Notons que si $P\in\Zr[x_1,x_2]$ est irréductible sur $\Zr$ (en particulier de contenu $1$) et est une équation de $V$, alors $\hat{h}(V)=\log M(P)$, où $M(P)$ est la mesure de Mahler de $P$. 

\medskip

Supposons l'in\'egalit\'e fausse pour une courbe $V$ de degr\'e $D$ définie sur $\Qr$, $\Qr$-irréductible qui ne soit pas de torsion. D'après un th\'eor\`eme de Zhang~\cite{Zh} on a $\me(V)\leq \frac{1}{D}\hat{h}(V)$, ainsi~:
\begin{equation}\label{contradiction}
\me(V)<\frac{1}{5^6D}\left(\frac{\log\log D'}{\log D'}\right)^3.
\end{equation}

\bigskip

\begin{center}
{\bf Choix des param\`etres et fonction auxiliaire}
\end{center}

On pose 
$$
\left\{\begin{disarray}{r @{\ :=\ } l}
T & \left[5 \frac{\log D'}{\log\log D'}\right]\\
L & DT^2\\
N & 5^4  \frac{(\log D')^2}{\log\log D'}.
\end{disarray}
\right.
$$

\bigskip

Nous utiliserons plusieurs fois l'inégalité suivante, valable pour tous réels $a,b,x>0$~:
\begin{equation}\label{MinExp}
\frac{x^a}{(\log x)^b}\geq \left(\frac{ea}{b}\right)^b.
\end{equation}

\bigskip

Notons que :
$$
T\geq [5e]\geq 13,\quad L\geq T^2\geq 169 \quad\text{et}\quad N\geq 5^4\frac{(\ln 16)^2}{\ln\ln 16}\geq 4000.
$$

\begin{fact}\label{5.1} On a 
$$
N/2\geq (\log D')^{1,99}
\quad\text{et}\quad
6,1\log(L+1)\leq T\log(N/2).
$$
\end{fact}

\bd
Pour la première inégalité, il suffit d'utiliser~(\ref{MinExp}) avec $(a,b)=(0.01,1)$. Pour la seconde, on a~:
$$
\frac{\log L}{T\log(N/2)} \leq  \frac{\log D'}{T\log(N/2)}+\frac{\log T}{T}\frac{2}{\log(N/2)}.
$$
De plus, comme $T\geq 13$, la premiére inégalité du fait nous donne $T\log(N/2)\geq 9,2\log D'$. Ainsi, puisque $N\geq 4000$~:
$$
\frac{\log L}{T\log(N/2)}\leq \frac{1}{9,2}+\frac{\log 13}{13}\frac{2}{\log(2000)}\leq\frac{1}{6,2}.
$$
Pour conclure, il suffit de remarquer que $\log (L+1)\leq 1,01\log L$, d'où le résultat.
\ed

\bigskip

En appliquant la proposition~\ref{Siegel2} à un ensemble fini suffisamment gros de points de $V$ de hauteur $\leq\theta$ où 
$$
\theta:=5^{-6} \left(\frac{\log\log D'}{\log D'}\right)^3,
$$
on trouve, par le même argument que celui utilisé dans le théorème~4.1 de~\cite{Am-Da2}, un polynôme non nul $F\in\Zr[\x]$, de contenu $1$ et de degré au plus $L$, qui s'annule en tout point de $V$ à un ordre $\geq T$ et tel que
$$
h(F) \leq l\Big\{(T+2)\log(L+1)+L\me(V)\Big\},
$$
où
$$
l:=\frac{{L+2 \choose 2}-{L-DT+2 \choose 2}}{{L-DT+2 \choose 2}}.
$$

\bigskip

\begin{center}
{\bf Extrapolation}
\end{center}

Soit $p$ un nombre premier dans $[N/2,N]$ et notons
$$
\begin{array}{lcl}
\varepsilon     & :=    & T\log p-h(F)-2\log(L+1) -pL\me(V)\\[2mm]
                & \geq  & T\log(N/2) -(l(T+2)+2)\log(L+1) - (N+l) L\me(V).
\end{array}
$$
La proposition~\ref{premiers} nous assure, via le même argument de densité que dans le lemme~4.2 de~\cite{Am-Da2}, que $F$ s'annule sur $[p]V$ si $\varepsilon>0\ $~; il suffit donc de montrer que notre choix de paramètres assure cette condition.

Majorons tout d'abord $l$ :
$$\label{MinHyper}
l \leq \frac{L+1}{L-DT+1}\times\frac{L+2}{L-DT+2}-1 \leq \left(\frac{L}{L-DT}\right)^2-1=\frac{2T-1}{(T-1)^2}
$$
en particulier, comme $T\geq 13$, nous obtenons :
$$
\begin{disarray}{lcl}
l(T+2)+2 	& \leq  & \frac{2T^2+3T-2}{(T-1)^2}+2\\[4mm]
		& = & \frac{2(T^2-2T+1)+7(T-1)+3}{(T-1)^2}+2\leq 5.
\end{disarray}
$$
Nous avons, puisque $l<1$ et $N>100$~:
$$
\varepsilon \geq T\log{N/2} - 5\log(L+1) - 1,01NL\me(V).
$$
Remarquons maintenant que, d'après le fait~\ref{5.1} nous avons $T\log(N/2)\geq 6,\!1\log(L+1)$, et d'après~(\ref{contradiction}) nous avons $NL\me(V)<\log D'\leq\log(L+1)$, ainsi $\varepsilon>0$.

\bigskip

\begin{center}
{\bf Conclusion}
\end{center}

Soit $\Lambda$ l'ensemble des nombres premiers dans $[N/2,N]$~; nous avons vu que, sous l'hypothèse~(\ref{contradiction}), $F$ s'annulait sur $[p]V$ pour tout premier $p\in\Lambda$. Comme $N\geq 4000$, la proposition~\ref{CardEcc} et le lemme~\ref{PiN} nous donnent :
\begin{equation}\label{finVe1}
L\geq \deg\Big(\bigcup_{p\in\Lambda} [p]V\Big)\geq \Big(0,4\frac{N}{\log N}-\frac{2}{\log 2}\log D\Big)D.
\end{equation}
Minorons le membre de droite~:
$$
\frac{\log\log D'}{\log N}\ \geq\ \frac{\log\log D'}{4\log 5+2\log\log D'}\ \geq \frac{3}{25},
$$
d'o\`u
$$
\frac{N}{\log N}\geq 3\left(\frac{5\log D'}{\log\log D'}\right)^2.
$$
En reportant ceci dans l'inégalité~(\ref{finVe1}) on obtient :
$$
\begin{disarray}{lcl}
L &\geq    & \left(    0,4\cdot 3 -\frac{2}{\log 2}\frac{(\log\log D')^2}{5^2\log D'}    \right) D\left(5\frac{\log D'}{\log\log D'}\right)^2\\[.4cm]
    &\geq    & \left(    1,2 -\frac{2}{\log 2}\frac{(\log\log 16)^2}{5^2\log 16}    \right) L\\[.4cm]
    & >     &  L
\end{disarray}
$$
d'o\`u une contradiction.


\subsection{Dans les extensions d'un corps cyclotomique}

F. Amoroso et U. Zannier donnent une minoration de la hauteur d'un nombre algébrique en fonction de son degré sur une extension abélienne de $\Qr$ ({\it c.f.}~\cite{Am-Za}) :

\begin{thm} Soient ${\mathbb K}$ un corps de nombres et $\Lr$ une extension abélienne de ${\mathbb K}$. Alors, pour tout $\gamma\in\Qbnt$, on a :
$$
h(\gamma) \geq \frac{C({\mathbb K})}{d} \left(\frac{\log \log 5d}{\log 2d}\right)^{13},
$$
où $d:=[\Lr(\gamma):\Lr]$ et $C({\mathbb K})$ est une constante dépendant uniquement de ${\mathbb K}$.
\end{thm}

Nous aurons besoin de ce résultat dans le cas particulier d'une extension cyclotomique, aussi nous nous proposons d'en montrer une version faible, mais explicite~: 

\begin{prop}\label{the:Dobrowolski2}
Soient $k=\Qr(\zeta_m)$ un corps cyclotomique et $\gamma\in\Qbnt$, alors :
$$
h(\gamma) \geq \frac{10^{-3}}{d} \left(\frac{\log \log D}{\log D}\right)^3,
$$
où $d:=[k(\gamma):k]$, $D:=[k(\gamma):\Qr]$ et $D\geq 2$.
\end{prop}

Notons $\varphi(m)$ l'indicatrice d'Euler de $m$~; nous pourrons supposer dans la suite $\varphi(m)=D/d\geq 3^4$. 
En effet, supposons le contraire, en particulier $[\Qr(\gamma):\Qr]\leq D=\varphi(m)d<3^4d$.\\ Si $[\Qr(\gamma):\Qr]\leq 15$, alors $h(\gamma)\geq 10^{-2}$, et si $16\leq[\Qr(\gamma):\Qr]<3^4d$, le théorème principal de~\cite{Vo} nous dit que~:
$$
h(\gamma) \geq \frac{1}{4[\Qr(\gamma):\Qr]} \left(\frac{\log \log [\Qr(\gamma):\Qr]}{\log [\Qr(\gamma):\Qr]}\right)^3.
$$
Par décroissance de la fonction $x\mapsto \frac{\log \log x}{\log x}$ sur $[16,+\infty]$ on en déduit
$$
h(\gamma) \geq \frac{1}{4D} \left(\frac{\log \log D}{\log D}\right)^3
\geq \frac{1}{324d}\left(\frac{\log \log D}{\log D}\right)^3.
$$
Notons toutefois qu'il est possible d'obtenir le même résultat sans~\cite{Vo}, avec toutefois une constante un peu plus petite : $2.10^{-4}$ au lieu $10^{-3}$ (la différence étant due à une moins bonne majoration de $\log (L+1)$ dans l'extrapolation).

\bigskip

\begin{center}
{\bf Choix des paramètres et fonction auxiliaire}
\end{center}

On pose
$$
\left\{
\begin{disarray}{lcl}
T       & = & \left[3\frac{\log D}{\log \log D}\right]\\[4mm]
L       & = & dT^2\\[4mm]
N       & = & 175\frac{(\log D)^2}{\log \log D}.
\end{disarray}
\right.
$$

\smallskip

\begin{fact}\label{factCyclo}
Nous avons~: $T\geq 8$, $L\geq 64$, $N\geq 175$, $T\leq 3D^{1/4}$ et 
$$
2,125\log (L+1) < 3,2\log D-\log \varphi(m).
$$
\end{fact}

\bd L'inégalité~(\ref{MinExp}) nous donne $T\geq [3e]=8$, (en particulier $L\geq 64$), et $N\geq 175\times 2e\geq 951$. Pour montrer $T\leq 3D^{1/4}$, comme $\log D\geq \log\varphi(m)\geq 4\log 3$, il suffit de remarquer que la fonction $x\mapsto \frac{x}{\log x}-e^{0.25x}$ est négative en $4\log 3$ et de dérivée négative sur $]1,+\infty[$.

Pour la dernière inégalité nous avons~:
$$
\begin{disarray}{lcl}
\log L  & = & \log D-\log\varphi(m)+2\log T\\[2mm]
        & \leq & 1,5\log D+2\log 3-\log\varphi(m)
\end{disarray}
$$
car $T\leq 3D^{1/4}$. Ainsi, puisque $\varphi(m)\geq 3^4$ il vient
$$
\log L\leq 1,5\log D-\frac{1}{2}\log \varphi(m).
$$
Comme $L\geq 64$ on a $\log (L+1)< 1,0038\log L$, ainsi
$$
2,125\log (L+1) < 3,2\log D-\log \varphi(m).
$$
\ed

Supposons par l'absurde~:
\begin{equation}\label{hcyclo}
h(\gamma) < \frac{10^{-3}}{d} \left(\frac{\log \log D}{\log D}\right)^3.
\end{equation}
D'après la proposition~\ref{Siegel2}, il existe un polynôme $F\in\Oc_k[{\mathrm x}]$ non nul, de degré au plus $L$, et nul en $\gamma$ à un ordre $\geq T$ tel que~:
$$
\begin{disarray}{lcl}
h(F)    & \leq & \frac{dT}{L+1-dT}\big(Lh(\gamma)+T\log(L+1)\big)+\log c_k\\[4mm]
        &\leq & \frac{1}{T-1}\big(Lh(\gamma)+T\log(L+1)\big)+\log \varphi(m).
\end{disarray}
$$
En effet, comme $\Delta_k|m^{\varphi(m)}$, on a 
$c_k=\sqrt{\frac{2}{\pi}\Delta_k^{1/\varphi(m)}}\leq \sqrt{\frac{2m}{\pi}}\leq \varphi(m)$.
Ainsi, comme $T\geq 8$ :
\begin{equation}\label{hgamma}
h(F)\leq Lh(\gamma)+1,125\log(L+1)+\log\varphi(m).
\end{equation}

\bigskip

\begin{center}
{\bf Extrapolation}
\end{center}

Soit maintenant $p\in[N/2,N]$ un nombre premier ne divisant pas $m$ (en particulier $p\nmid\Delta_k$) et soit $\nu$ est une place ne divisant pas $p$. Notons $T_1$ l'ordre d'annulation de $F^{\phi_p}$ en $\gamma^p$~; nous devons donc montrer $T_1>0$. D'après la proposition~\ref{premiers}, on a $T_1\Big(\log (L+1)+\log p\Big)>\varepsilon$ où 
$\varepsilon := T\log p-h(F)-pLh(\gamma)-\log(L+1)$, il nous suffit ainsi de voir que $\varepsilon$ est strictement positif. Par hypothèse sur $F$ on a~: 
$$
\varepsilon \geq T\log p-Lh(\gamma)\left(1+p\right)-2,125\log(L+1)-\log\varphi(m)
$$
soit, d'après le fait~\ref{factCyclo}~:
$$
\varepsilon > T\log p-Lh(\gamma)(1+p)-3,2\log D.
$$
Nous allons voir que $F^{\phi_p}$ s'annule en $\gamma^p$, en montrant que $\varepsilon>0$.
Par hypothèse $p\geq N/2\geq \big(\log D\big)^{1,99}$ (d'après~(\ref{MinExp}) avec $x=\log D$, $a=0,01$ et $b=1$) et $T\geq 8$, ainsi :
$$
\begin{disarray}{rcl}
\varepsilon     & > & \Big(\frac{8}{9}(T+1)1,99\log\log D -3,2\log D\Big)-Lh(\gamma)(1+p)\\[4mm]
        & > & 2\log D-Lh(\gamma)(1+N)\\[4mm]
        & > & 2\log D-176\cdot 3^2d\frac{(\log D)^4}{(\log \log D)^3}\times h(\gamma).
\end{disarray}
$$
Donc, d'après~(\ref{hcyclo}), on a $\varepsilon>0$, en particulier $T_1>0$ et $F^{\phi_p}(\gamma^p)=0$.

\bigskip

\bigskip

\begin{center}
{\bf Conclusion}
\end{center}

\noindent  Remarquons que si $F^{\phi_p}(\gamma^p)=0$ et si $\tau\in$ Gal$(\Qb/k)$ prolonge $\phi_p^{-1}$,  alors $F(\tau(\gamma^p))=0$. Notons $\Sigma$ l'ensemble des $\tau(\gamma^p)$, où 
\begin{itemize}
\item $p$ parcourt l'ensemble des premiers $p$ de $[N/2,N]$ ne divisant pas $m$,
\item $\tau\in$ Gal$(\Qb/k)$ prolonge $\phi_p^{-1}$.
\end{itemize}
Sous l'hypothèse~(\ref{hgamma}) sur la hauteur de $\gamma$, nous avons vu que $F$ s'annule sur $\Sigma$. 
Pour arriver à une contradiction, nous allons montrer que $\card(\Sigma)>\deg(F)$.

$\bullet$ Soient $p_1<\dots<p_s$ les diviseurs premiers de $m$ dans $\{N/2,\ldots,N\}$, nous avons :
$$
\varphi(m)\geq (p_1-1)(p_2-1)\cdots (p_s-1)\geq (N/2-1)^s,
$$
en particulier
$$
s\leq\frac{\log(\varphi(m))}{\log (N/2-1)}\leq \frac{\log D}{6}.
$$

$\bullet$ Remarquons maintenant que si $p$ est un nombre premier tel que $\Qr(\gamma^p)=\Qr(\gamma)$, alors $k(\gamma^p)=k(\gamma)$ c'est-à-dire $[k(\gamma^p):k]=[k(\gamma):k]$, auquel cas les $k$-automorphismes de $\Qb$ prolongeant $\phi_p^{-1}$ sont au nombre de $d=[k(\gamma):k]$. 

Ainsi, d'après la proposition~\ref{CardEcc} page~\pageref{CardEcc}, le nombre de premiers $p$ tels que $[k(\gamma^p):k]<[k(\gamma):k]$ ou tel que $\gamma^p$ soit égal $\sigma(\gamma^p)$ pour un certain $\sigma\in$ Gal$(\Qb/k)$, $\sigma\neq id$, est inférieur à :
$$
\frac{\log([\Qr(\gamma):\Qr])}{\log 2}\leq \frac{\log D}{\log 2}.
$$
De plus, comme $\gamma$ n'est pas racine de l'unité, $\gamma^p$ et $\gamma^{p'}$ ne sont pas conjugués si $p\neq p'$, d'où :
$$
\card(\Sigma)\geq d\times \left(\pi(N)-\pi(N/2)-\left(\frac{1}{\log 2}+\frac{1}{6}\right)\log D\right)
$$
soit, d'après le lemme~\ref{PiN} page~\pageref{PiN} et~(\ref{MinExp}) avec $x=\log D$, $a=1$ et $b=2$ :
$$
\begin{disarray}{lcl}
\card(\Sigma)   & \geq & d\left( 0,\!4\frac{N}{\log N}-1,\!7\log D\right)\\[4mm]%
                & \geq & d\left( 0,\!4\frac{N}{\log N}-1,\!7\left(\frac{2}{e}\right)^2\left(\frac{\log D}{\log\log D}\right)^2\right).
\end{disarray}
$$
De plus, $\log N\leq (2+\log 175)\log\log D$, donc :
$$
\card(\Sigma)\geq d\left( \frac{0,4\cdot 175}{2+\log 175}-1\right)\left(\frac{\log D}{\log\log D}\right)^2,
$$
d'où
$$
\card(\Sigma)> 3^2d\left(\frac{\log D}{\log\log D}\right)^2.
$$
Ainsi $\card(\Sigma)>L\geq \deg(F)\ $ ; en particulier $F$ ne peut pas s'annuler sur $\Sigma$ tout entier, contradiction avec l'hypothèse~(\ref{hcyclo}). \hfill $\blacksquare$

\bigskip
\section{Démonstration du théorème principal}\label{DemoMthm}

Si $V_{\alphab}$ désigne la variété de dimension zéro définie sur $\Qr$ par un point $\alphab$ de $\Gm^n$, c'est-à-dire $\{\sigma(\alphab)\ |\ \sigma\in\text{Gal}(\Qb/\Qr)\}$, l'inégalité~(\ref{deg_obs}) nous dit
$$
\io_{\Qr}(\alphab)\leq n(\deg(V_{\alphab}))^{1/n}.
$$
De façon générale, si $V$ est une variété définie sur $\Qr$ et $\Qr$ irréductible contenant $\alphab$, il découle immédiatement d'un résultat de {\sc M. Chardin} ({\it cf.}~\cite{Ch}, 
corollaire~2, chapitre~1, page~310 et exemple~1, page~311) l'inégalité :
\begin{equation}
\label{deg_obsCH}
\io_{\Qr}(\alphab)\leq n\deg(V)^{1/\codim(V)}.
\end{equation}
Afin de pouvoir conclure leur démonstration, \cite{Am-Da4} considèrent des variétés de différentes dimensions contenant un translaté de la variété qu'ils étudient, aussi ont-ils été amenés à introduire {\it l'indice d'obstruction généralisé de poids $T$}~:
$$
\io(T;\alphab):=\min\{(T\deg(W))^{1/\codim(W)}\}\ ,
$$
où $T$ est un réel $>0$ et $W$ parcourt l'ensemble des variétés définies sur $\Qb$ et $\Qb$ irréductibles contenant $\alphab$. Notons en particulier qu'aucune hypothèse sur le corps de définition de $V$ n'est faite ici. Nous utiliserons cet {\it indice d'obstruction généralisé} un peu modifié, en gros~:
$$
\min\left\{\ 2\io T,\ (TD)^{1/2}\right\}\ ,
$$
où $D$ est le degré sur $\Qr$ d'un point $\alphab$ et le $\io$ le degré d'une courbe $V$ définie sur $\Qr$ fixée contenant $\alphab$ (voir dans le choix des paramètres ci-dessous). Celui-ci dans notre cas n'est pas nécessaire pour retrouver la minoration du théorème~\ref{Mthm} à une puissance du log près, néanmoins il permet de gagner non seulement sur la constante, mais surtout dans le terme d'erreur (sur la puissance du $\log$).

\subsection{Choix des paramètres et fonction auxiliaire}

Notons $D$ le degré de $\Qr(\alphab)$ sur $\Qr$ et posons~:
$$
\left\{
\begin{disarray}{lcl}
T & = & \left[9\left(\frac{\log  \io'}{\log\log  \io'}\right)^2\right]\\[5mm]   
L & = & \min\left\{\ 2\io T^2,\ \left[(TD)^{1/2}(T+1)\right]\right\}.
\end{disarray}
\right.
$$

Notons $c_1:=3,\!7\cdot 10^4$, $c_2:=2,\!05\cdot 10^{9}$ et considérons les réels $N_1,N_2$ suivants :
$$
\left\{
\begin{disarray}{lcl}
N_1 & := & c_1 \frac{\big(\log\io'\big)^2}{\log\log\io'}\\[4mm]
N_2 & := & c_2 \frac{\big(\log\io'\big)^{8}}{\big(\log\log\io'\big)^6}. 
\end{disarray}
\right.
$$

\begin{fact}\label{factThmPrinc}
Nous avons~:
\begin{enumerate}
\item $T\geq [9e^2]=66$ et $N_1^2\leq N_2$.
\item $\log (L+1)\leq 4,\!3\log\io.$
\item $\log (N_1/2)\geq 1,\!999\log\log \io'$ et $\log(N_2/2)\geq 7,\!92\log\log \io'$.
\item Les inégalités suivantes nous permettrons de majorer le cardinal d'ensembles de premiers exceptionnels~:
$$
\frac{2}{\log 2}\log L\leq 0,\!01\frac{N_1}{\log N_1}\quad\text{et}\quad
\frac{2}{\log 2}\log (N_1L^2)\leq 0,\!01\frac{N_2}{\log N_2}.
$$
\end{enumerate}
\end{fact}

\bd
\begin{enumerate}
\item Pour la première inégalité on utilise~(\ref{MinExp}) avec $a=b=2$. La seconde découle immédiatement du choix des constantes.
\item Comme $T\geq 66$ nous avons $L+1\leq 2\io'T^2+1\leq 1,0002\cdot 2\io'T^2$ or~: 
$$
\begin{array}{lcl}
\log L	& \leq & \log 2+\log \io'+2\log (9/\log\log 16)+4\log\log \io'\\[3mm]
	& \leq & \Big(1+(\log 2+2\log 8,9+4\log\log 16)/\log 16\Big)\log \io'\\[3mm]
	& \leq & 4,299\log\io'
\end{array}
$$
d'où le résultat.
\item Il suffit de remarquer que, d'après~(\ref{MinExp}) avec $(a,b)=(0.001,1)$ puis $(0.08,6)$ nous avons~:
$$
\frac{N_1}{2}\geq\frac{c_1}{2}\cdot 0,\!001e(\log\io')^{1,99}\quad\text{et}\quad
\frac{N_1}{2}\geq\frac{c_1}{2}\left(\frac{0,\!08e}{6}\right)^6(\log\io')^{7,92}.
$$
\item Comme $\log N_1\leq (2+\log c_1)\log\log \io'$ on a 
$$
0,\!01\frac{N_1}{\log N_1}\geq \frac{c_1}{100(2+\log c_1)} \left(\frac{\log\io'}{\log\log\io'}\right)^2 
\geq \frac{2\cdot 4,3\log\io'}{\log 2}\frac{\log\io'}{\big(\log\log\io'\big)^2}$$
or $4,3\log\io'\geq \log L$ et d'après~(\ref{MinExp}), $\frac{\log\io'}{\big(\log\log\io'\big)^2}\geq \frac{e^2}{4}\geq 1$, d'où l'inégalité voulue.

Enfin nous avons, puisque $N_2\geq N_1^2$~:
$$
0,\!01\frac{N_2}{\log N_2}\geq 0,\!01N_1\frac{N_1}{2\log N_1}\geq 100\frac{N_1}{\log N_1}
$$
de plus
$$
\frac{2}{\log 2}\log (N_1L^2)= \frac{2}{\log 2}\left(\log N_1+2\log L\right)\leq \frac{2}{\log 2}\log N_1+0,\!02\frac{N_1}{\log N_1}.
$$
\end{enumerate}
\ed

\medskip

Le corollaire~\ref{f_auxCor} nous donne un polynôme $F\in E(\{\alphab\},L,T)\cap \Zr[\x]$ vérifiant : 
$$
h(F) \leq \frac{1}{T-1}\big( (T+1)\log (L+1)+Lh(\alphab)\big)
$$
Pour $j=1,2$ posons $\Pc_j:=\{p\in [N_j/2,N_j]\ \text{premier}\}\cup\{1\}$ et notons
$$
T_1:=\min_{p\in\Pc_1} \text{ord}_{\alphab^p}(F),
$$
en particulier, comme $1\in\Pc_1$ nous avons $T_1\leq T$. Nous allons montrer le théorème~\ref{Mthm} par l'absurde, aussi nous supposerons dans la suite $\alphab$ de hauteur petite, plus précisément :

\noindent\rule{\textwidth}{.2mm}

\vspace{-4mm}

\begin{equation}\label{hPetit}
h(\alphab)\leq\frac{T_1\log N_2}{10N_1N_2L}\leq\frac{(T+1)}{10L}.
\end{equation}

\vspace{-2mm}

\noindent\rule{\textwidth}{.2mm}

On a alors :
$$
h(F) \leq \frac{T+1}{T-1}\left( 1+\frac{1}{10\log (L+1)} \right)\log (L+1)
$$
ainsi, comme $T\geq 66$ et $L+1\geq T^{3/2}\geq 500$ on obtient
$$
h(F)\leq 1,05\log (L+1).
$$

\subsection{Extrapolation}\label{MthmExrtap}

\begin{prop}
Sous l'hypothèse~(\ref{hPetit}) sur la hauteur de $\alphab$, nous avons, pour tout $(p,q)$ dans $\Pc_1\times\Pc_2$,
$$
F(\alphab^{pq})=0.
$$
\end{prop}

\bd
Soit $(p,q)\in\Pc_1\times\Pc_2$, puisque $T_1\leq T$ et $N_1\geq 100$ nous avons d'après~(\ref{hPetit})
$$
Lh(\alphab)\leq\frac{T_1}{10N_1}\frac{\log N_2}{N_2}\leq\frac{T}{1000}\frac{\log N_2}{N_2}.
$$
Par décroissance de la fonction $x\mapsto \log(x)/x$, comme $p\leq N_1\leq N_2/2\leq q\leq N_2$ il découle
\begin{equation}\label{hPetit2}
pLh(\alphab)\leq 0,\!001 T\log p\quad\text{et}\quad pqLh(\alphab)\leq 0,\!1 T_1\log q.
\end{equation}
$\bullet$ Notons $T_{1,p}$ l'ordre d'annulation de $F$ en $\alphab^p$, comme $h(F)\leq 1,\!05\log(L+1)$, la proposition~\ref{premiers} nous donne :
$$
T_{1,p} \big( \log (L+1)+\log p\big)>0,\!999 T\log p-3,\!05\log (L+1).
$$
Deux cas apparaissent :

\underline{si $L+1\leq p$ :}
\begin{equation}\label{Lpetit}
2T_{1,p}\log p>(0,\!999T-3,\!05)\log p,
\end{equation}      
donc, comme $T\geq 66$, on obtient $T_{1,p}\geq 32$.      

\underline{si $L+1>p$ :}
\begin{equation}\label{ordT_1b}
(2T_{1,p}+3,05)\log (L+1) > 0,\!999T\log p \geq 0,\!999T\log (N_1/2)
\end{equation}
Comme $T\geq 66$, on a, d'après le point 2 du fait~\ref{factThmPrinc}~:
$$
\begin{disarray}{rcl}
2T_{1,p}+3,\!05      & > & \frac{0,\!999T}{T+1}(T+1)\frac{\log (N_1/2)}{4,\!3\log  \io'}\\[4mm]
                & > &  \frac{0,\!999\cdot 66}{67}\cdot 9\cdot \frac{\log
  (N_1/2)}{4,\!3}  \frac{\log\io'}{(\log\log\io')^2}
\end{disarray}
$$
Soit, en utilisant les inégalités $N_1/2\geq \frac{c_1\log 16}{2(\log\log 16)^2}$ et $\frac{\log\io'}{(\log\log\io')^2}\geq \frac{e^2}{4}$ (via (\ref{MinExp}))~:
$$
2T_{1,p}+3,05 > 45,\!05
$$
car $c_1=3,7.10^4$. Ainsi dans les deux cas on a $T_{1,p}\geq 22$.

\medskip

$\bullet$ Notons maintenant $T_{2,pq}$ l'ordre d'annulation de $F$ en $\alphab^{pq}$\ ; nous avons d'après~(\ref{hPetit2}) $pLh(\alphab^q)=pqLh(\alphab)\leq 0,\!1 T_{1,p}\log q$. Comme $h(F)\leq 1,\!05\log(L+1)$, de nouveau la proposition~\ref{premiers} nous donne :
$$
T_{2,pq} \big( \log (L+1)+\log q\big)>0,\!9T_{1,p}\log q-3,\!05\log (L+1).
$$
Il nous faut ici montrer que le membre de droite de cette inégalité est $>0$. Si $L+1\leq p$, c'est évident, car $T_{1,p}\geq 22$ et $q\geq p$. Supposons donc $L+1>p\ $ ; comme $T_{1,p}\geq 22$ nous avons :
$$
0,\!9T_{1,p}> 0,\!42(2T_1+3,\!05),
$$
ainsi d'après~(\ref{ordT_1b})
$$
0,\!9T_{1,p}\log (N_2/2) > 0,\!42\frac{0,\!999T\log (N_1/2)}{\log (L+1)}\log (N_2/2).
$$
D'où, d'après les points 1,2 et 3 du fait~\ref{factThmPrinc}~:
\begin{equation}\label{Lgrand}
\begin{array}{lcl}
0,\!9T_{1,p}\log (N_2/2)	& > & 0,\!42\cdot 0,\!999 \cdot \frac{66}{67}\cdot \frac{1,\!999}{4,\!3}\cdot 7,\!92\cdot 9\log  \io'\\[2mm]
        & > & 13,\!5\log \io'.
\end{array}
\end{equation}
Ainsi $0,\!9T_{1,p}\log q> 3,\!05\cdot 4,\!3 \log \io'\geq 3,\!05\log(L+1)$.

\ed

\subsection{Conclusion}\label{sub5.3}

Notons $X_1$ la variété définie par $F$ et posons :
$$
\left\{
\begin{disarray}{l@{\ :=\ } l}
X_2 & \bigcap_{p\in\Pc_1} [p]^{-1}X_1,\\[4mm]
X_3 & \bigcap_{(p,q)\in\Pc_1\times\Pc_2} [pq]^{-1}X_1.
\end{disarray}
\right.
$$
Notons que, puisque $\Pc_1$ et $\Pc_2$ contiennent $1$, nous avons les inclusions suivantes~:
$$
X_3\subset X_2\subset X_1.
$$
Nous travaillons ici avec $\alphab$, aussi nous ne considérerons que les composantes de ces variétés rencontrant une puissance de $\alphab$~; plus précisément posons :
\begin{itemize}
\item $Y_1$ l'union des composantes $\Qr$-irréductibles de $X_1$ contenant $\alphab^{pq}$ pour au moins un $(p,q)$ dans $\Pc_1\times\Pc_2$
\item $Y_2$ l'union des composantes $\Qr$-irréductibles de $X_2$ contenant $\alphab^{q}$ pour au moins un $q\in\Pc_2$
\item $Y_3$ l'union des composantes $\Qr$-irréductibles de $X_3$ contenant $\alphab$.
\end{itemize}

On a les inclusions suivantes :
$$
\alphab\in Y_3\subset Y_2\subset Y_1
$$
En particulier, deux de ces trois variétés ont même dimension ce qui nous permettra de comparer leurs degrés ou leur hauteurs normalisées.

\subsubsection{Cas où $Y_2$ et $Y_3$ sont de dimension $0$}\label{CasDim0}

Soit $Z$ une composante $\Qr$-irréductible de $Y_3$\ ; comme $Z$ rencontre $\alphab$ on a~:
$$
Z=\bigcup_{\sigma\in\text{Gal}(\Qb/\Qr)}\sigma(\alphab).
$$
En particulier $\deg(Z)=D$. De l'inclusion 
$$
\bigcup_{q\in\Pc_2} [q]Z \subseteq Y_2\ ,
$$
on obtient une première inégalité :
\begin{equation}\label{dim0}
\deg\Big(\bigcup_{q\in\Pc_2} [q]Z\Big) \leq \deg Y_2.
\end{equation}

$\pg$ Soient $F_1,\dots,F_r$ les facteurs $\Qr$-irréductibles de $F$. Les composantes $\Qr$-irréductibles de $X_2$ de dimension $1$ sont les $Z(F_j)$, où :
$$
F_j\mid \gcd \big(\{F(\x^p),\ p\in\Pc_1\}\big).
$$
Quitte à les  réordonner, on peut supposer que $1,\ldots,l$ sont les indices $i$ pour lesquels $F_i$ ne divise pas $\gcd \big(\{F(\x^p),\ p\in\Pc_1\}\big)$. En particulier, comme $Y_2$ est de dimension zéro par hypothèse, si $j\in\{l+1,\ldots,r\}$, alors $F_j(\alphab^q)\neq 0$ pour tout $q\in\Pc_2$. Choisissons maintenant un polynôme $G$ de la forme 
$$
G(\x)=\sum_{p\in \Pc_1\setminus\{1\}} \lambda_p F(\x^p) \qquad \lambda_p\in\Qr
$$
tel que $G$ ne soit pas un diviseur de zéro de $\Qr[\x]/(\tilde{F})$. Un tel polynôme existe bien, il suffit en effet de remarquer que pour tout $j$ dans $\{1,\dots,l\}$, le sous-espace vectoriel 
$$
\left\{\lambdab\in\Qr^{\card(\Pc_1)-1}\ \Big|\ \sum_{p\in\Pc_1\setminus\{1\}} \lambda_p F(\x^p)\in (F_j)\right\}
$$
est propre. Comme $Y_2\subseteq Z(\tilde{F})\cap Z(G)$ et ce dernier est de dimension $0$, le théorème de Bézout nous donne~:
$$
\deg Y_2\leq \deg(F)\deg(G)\leq  N_1L^2\leq N_1 DT(T+1)^2.
$$

$\pg$ Considérons maintenant le membre de gauche de~(\ref{dim0}). Comme $N_2\geq 5000$, la proposition~\ref{CardEcc} et le lemme~\ref{PiN} nous donnent :
$$
\left(0,\!41\frac{N_2}{\log N_2}-\frac{2}{\log 2}\log D\right)D\leq \deg\Big(\bigcup_{q\in\Pc_2} [q]Z\Big).
$$
De plus, comme $Z\subset Y_2$, on a $D\leq\deg(Y_2)\leq N_1L^2$, soit, d'après le point 4 du fait~\ref{factThmPrinc}~:
$$
\frac{2}{\log 2}\log D\leq \frac{2}{\log 2}\log (N_1L^2)\leq 0,\!01\frac{N_2}{\log N_2}.
$$
En reportant tout ceci dans~(\ref{dim0}) on en déduit
$$
\frac{0,\!4N_2}{\log N_2}\leq N_1 T(T+1)^2.
$$
D'où, en utilisant les inégalités $\log N_2\leq (8+\log c_2)\log\log\io'$ et $T\geq 66$~:
$$
\frac{0,\!4c_2}{8+\log c_2}\leq c_1 9^3 \left(\frac{67}{66}\right)^2,
$$
contradiction, car $c_1=3,\!7\cdot 10^4$ et $c_2=2,\!05\cdot 10^9$.

\bigskip

\subsubsection{Cas où $Y_1$ et $Y_2$ sont de dimension $1$}\label{CasDim1}

Soit $Z$ une composante $\Qr$-irréductible de $Y_2$ de dimension $1$, et soit $q\in\Pc_2$ tel que $\alphab^q\in Z$.

\medskip

\noindent $\bullet$ Supposons dans un premier temps que $Z$ soit de torsion. Si $B$ désigne $[q]^{-1}Z$, alors $\alphab\in B$ et $B$ est de torsion. Comme $Z$ et $Y_1$ sont de même dimension, on a :
$$
\deg(B)\leq N_2\deg(Z)\leq N_2\deg(Y_1)\leq N_2L\leq 2c_2 9^2\io \big(\log \io'\big)^{12}.
$$
De plus, $V$ étant irréductible et non de torsion, $V$ et $B$ n'ont pas de composante commune, le théorème de Bézout nous donne~:
$$
D\leq \deg(V)\cdot \deg(B)= \io\deg(B),
$$
où $D$ est le degré de $\alphab$ sur $\Qr$. Ainsi, comme $\io'\geq 16$~:
$$
D\deg(B)\leq (2c_2 9^2)^2\io^3\big(\log \io'\big)^{24}\leq  \io'^{23}\big(\log \io'\big)^{24}\leq  \io'^{47}.
$$
Le lemme~\ref{torsion} ci-dessous nous dit alors :
$$
h(\alphab)\geq\frac{10^{-9}}{\io}\left(\frac{\log\log \io'}{\log \io'}\right)^{3}.
$$
ce qui nous donne bien le théorème~\ref{Mthm}.

\bigskip

\begin{lem}\label{torsion}
Soit $V$ une courbe de $\Gm^2$ définie sur $\Qr$ et $\Qr$-irréductible de degré $\io$ qui ne soit pas de torsion, et soit $\alphab\in V\setminus(\Gm^2)_{\text{tors}}$. S'il existe une courbe $B$ de torsion définie et irréductible sur $\Qr$ contenant $\alphab$, alors~:
$$
h(\alphab)\geq \frac{5.10^{-4}}{\io}\left(\frac{\log(D\deg(B))}{\log\log (D\deg(B))}\right)^{-3}.
$$
\end{lem}

\bd  
Il existe un sous-groupe algébrique $H$ de $\Gm^2$ et un $\thetab\in(\Gm^2)_{\text{tors}}$ tels que~:
$$
B=\bigcup_{\sigma\in\text{Gal}(\Qb/\Qr)} \sigma(\thetab) H.
$$
Soient $a,b\in\Zr$ tels que~:
$$
H=\big\{ (x,y)\in\Gm^2\ |\ x^a y^b=1 \big\}.
$$
Comme $\alphab\in B$, il existe $\eta\in(\Gm)_{\text{tors}}$ tel que $\alpha_1^a\alpha_2^b=\eta$. Soit $\gamma$ une racine $b$-ième de $\alpha_1$ ($\alphab$ n'étant pas de
torsion, $\gamma\not\in(\Gm)_{\text{tors}}$), on a $\alpha_2^b =
\eta\gamma^{-ab}$. En particulier, il existe
$\eta'\in(\Gm)_{\text{tors}}$ tel que $\alpha_2=\eta'\gamma^{-a}$. Posons $M:=\max\{|a|,|b|\}$, on a :
$$
\begin{array}{lcl}
h(\alphab) & \geq & \max\{h(\alpha_1),h(\alpha_2)\}\\
             & \geq & \max\{h(\gamma^b),h(\eta'\gamma^{-a})\}\\
             & \geq & M\cdot h(\gamma)  
\end{array}
$$
Considérons
$$
g(t):=t^\lambda G(t^b,\eta t^{-a})\in\Qr(\eta)[t],
$$
où $G\in\Qr[\x]$ est une équation de de $V$ et $\lambda\in\Zr$ est choisi le plus petit possible. En particulier $G$ est nul en $\alphab$ de degré $\io$ et {\em a fortiori} on a $g(\gamma)=0$. Notons que, comme $V$ n'est pas de torsion, le polynôme $g$ est non nul.

Notons $D_\gamma:=[\Qr(\eta,\gamma):\Qr]$ et $d_\gamma:=[\Qr(\eta,\gamma):\Qr(\eta)]$~; l'extension $\Qr(\eta)/\Qr$ étant cyclotomique, la proposition~\ref{the:Dobrowolski2} nous dit que~:
$$
h(\gamma)\geq \frac{10^{-3}}{d_\gamma}\left(\frac{\log(D_\gamma)}{\log\log (D_\gamma)}\right)^{-3},
$$
or $d_\gamma\leq \deg(g) \leq 2 \deg(G) M = 2\io M$ et $D_\gamma\leq DM$,
d'où :
$$
h(\gamma)\geq \frac{5.10^{-4}}{\io M}\left(\frac{\log(MD)}{\log\log (MD)}\right)^{-3},
$$
ainsi :
$$
\begin{disarray}{lcl}
h(\alphab) & \geq &  M\cdot h(\gamma)\\[3mm]
	& \geq &  \frac{5\cdot 10^{-4}}{\io}\left(\frac{\log(MD)}{\log\log (MD)}\right)^{-3}.\\
\end{disarray}
$$
Pour finir, il suffit de remarquer que, comme $H$ et $B$ ont la même dimension, on a 
$$
M\leq\deg(H)\leq\deg(B).
$$
\ed

\bigskip

\noindent $\bullet$ Nous supposerons dans la suite que $Z$ n'est pas de torsion. 
Nous avons l'inclusion 
$$
\bigcup_{p\in\Pc_1} [p]Z \subseteq Y_1.
$$
Comme les variétés $Z$ et $Y_1$ sont de même dimension, on en déduit~:
$$
\hat{h}(Y_1)\geq \hat{h}\Big(\bigcup_{p\in\Pc_1} [p]Z\Big)
$$
Notons $W_1,\dots,W_l$ les composantes géométriquement irréductibles de $Z$. Comme $Z$ n'est pas de torsion, le lemme~2.3 de~\cite{Am-Da} nous dit que,  si $(p,i)$ et $(p',j)$ sont deux couples distincts d'éléments de  $\left(\Pc_1\setminus Ecc(Z)\right)\times\{1,\dots,l\}$, alors les sous-variétés $[p]W_i$ et $[p']W_j$ sont distinctes~; ainsi
\begin{equation}\label{htnpZ}
\hat{h}(Y_1)\geq \sum_{p\in\Pc_1\atop p\not\in Ecc(Z)}\sum_{i=1}^l \hat{h}([p]W_i).
\end{equation}
Si $W$ désigne une composante géométriquement irréductible de $Z$, il nous faut donc majorer le cardinal de $Ecc(Z)$ et évaluer $\hat{h}([p]W)$ en fonction de $\hat{h}(W)$. Rappelons que le stabilisateur de $W$ est par définition~:
$$
G_W:=\{\y\in\Gm^n\ |\ \y\cdot W=W\}=\bigcap_{y\in W} y^{-1}W,
$$
en particulier $\dim(G_W)\leq\dim(W)=1$. Notons ici que les premiers divisant le cardinal de $G_W/G_{W}^0$ (quotient de $G_W$ par sa composante neutre $G_W^0$) sont dans $Ecc(W)$\footnote{cardinal qui est indépendant du choix de la composante $W$.}. On sait de plus, d'après la proposition 2.1 de~\cite{Da-Ph} que~:
$$
\hat{h}([p]W)=\frac{p^{\dim(W)+1}}{|\ker [p]\cap G_W|}\hat{h}(W),
$$
et $|\ker [p]\cap G_W|=p^{\dim G_W}|\ker [p]\cap G_W/G_W^0|\leq p\cdot|\ker [p]\cap G_W/G_W^0|$. En particulier si $p\not\in Ecc(Z)$, auquel cas $p$ ne divise pas $|G_W/G_W^0|$, on a~:
$$
\hat{h}([p]W)\geq p\cdot\hat{h}(W).
$$
La proposition~\ref{CardEcc} et le point 4 du fait~\ref{factThmPrinc} nous donnent de plus
$$
\card(Ecc(Z)\cap\Pc_1)\leq\frac{2\log\deg(Z)}{\log 2}\leq \frac{2\log L}{\log 2}\leq 0,\!01\frac{N_1}{\log N_1}.
$$ 
Ainsi, en reportant ceci dans~(\ref{htnpZ})~:
$$
\hat{h}(Y_1) \geq \sum_{p\in\Pc_1\atop p\not\in Ecc(Z)} p\cdot\hat{h}(Z)
    \geq \left(\pi(N_1)-\pi(N_1/2)-0,\!01\frac{N_1}{\log N_1}\right)\frac{N_1}{2} \cdot\hat{h}(Z).
$$
Comme $N_1\geq 5000$, nous déduisons du lemme~\ref{PiN}~:
\begin{equation}\label{HNinter}
\hat{h}(Y_1)\geq 0,\!2\frac{N_1^2}{\log N_1}\cdot\hat{h}(Z).
\end{equation}
Comme $\dim Z=\dim Y_1=\dim X_1=1$ et $Z\subset Y_1\subset X_1$ on a $\deg(Z)\leq \deg Y_1\leq L$. La variété  $Z$ n'étant pas de torsion, la proposition~\ref{dim2} nous dit :
$$
\hat{h}(Z)\geq 5^{-6} \left(\frac{\log\log L}{\log L}\right)^3,
$$
de plus, l'inégalité de Landau  $h(F)+\log \deg(F)\geq \log M(F)$ nous donne :
$$
2,\!05\log L\geq h(F)+\log \deg(F)\geq\log M(F)=\hat{h}(X_1)\geq \hat{h}(Y_1).
$$
En reportant tout cela dans~(\ref{HNinter}) on obtient alors
$$
2,\!05\log L\geq 5^{-7} \frac{N_1^2}{\log N_1} \left(\frac{\log\log L}{\log L}\right)^3.
$$
Remarquons maintenant que $\log N_1\leq (2+\log c_1)\log\log\io'$ et, d'après le point 2 du fait~\ref{factThmPrinc}, que $\log L\leq 4,\!3\log \io'$, soit
$$
\frac{N_1^2}{\log N_1} \geq \frac{c_1^2}{2+\log c_1}\frac{\big(\log\io'\big)^4}{\big(\log\log \io'\big)^3}> 2,\!05\cdot 5^7\cdot 4,\!3^4\cdot \frac{\big(\log  \io'\big)^4}{\big(\log\log \io'\big)^3}
$$
car $c_1=3,\!7\cdot 10^4$, contradiction. 

\bigskip

\subsubsection{Conclusion de la démonstration du théorème~\ref{Mthm}}

L'hypothèse~(\ref{hPetit}) est donc fausse, ainsi :
$$
h(\alphab) \geq \frac{T_1\log N_2}{10N_1N_2L}.
$$
\begin{fact}\label{fact4.3}
On a 
$$
T_1\log N_2\geq 15\log  \io'
$$
\end{fact}

\bd
Soit $(p,q)\in\Pc_1\times\Pc_2$, tel que l'ordre d'annulation $T_{1,p}$ de $F$ en $\alphab^p$ vérifie $T_{1,p}= T_1$. 
Si $L+1>p$, c'est exactement l'inégalité~(\ref{Lgrand}). Supposons donc $L+1\leq p$, comme $N_2\geq N_1^2$ et $N_1\geq p\geq N_1/2$, l'inégalité~(\ref{Lpetit}) donne :
$$
T_1\log N_2\geq 2T_1\log N_1 > (0,\!999T-3,\!05)\log(N_1/2)
$$
soit, d'après le point 3 du fait~\ref{factThmPrinc}~:
$$
\begin{array}{lcl}
T_1\log N_2	& > & 1,\!999(0,\!999T-3,\!05)\log\log  \io'\\
		& > &  (T+1)\log\log \io'\\[1mm]
		& > &  9\frac{(\log  \io')^2}{\log\log \io'}\\[1mm]
                        & > & 9e\log \io'
\end{array}
$$
\ed
Par définition, $L\leq 2\io T^2$, le Fait~\ref{fact4.3} montre alors que
$$
h(\alphab) \geq \frac{1,\!5}{2\cdot 9^2c_1c_2}\frac{1}{\io}\frac{(\log\log  \io')^{11}}{(\log  \io')^{13}}
$$
ainsi,
$$
h(\alphab) \geq \frac{1,\!2\cdot 10^{-16}}{\io} \frac{\big(\log\log  \io'\big)^{11}}{\big(\log  \io'\big)^{13}}
$$
car $c_1=3,\!7.10^4$ et $c_2=2,\!05.10^{9}$.


\section{Démonstration du corollaire~\ref{Mcor}}

Soient $\alphab:=(\alpha_1,\alpha_2)$ un point à coordonnées multiplicativement indépendantes. On peut supposer sans perte de généralité $h(\alpha_1)\leq h(\alpha_2)\leq 1$ et $D:=[\Qr(\alphab):\Qr]\geq 2$. Soient $A\in\Nr^*$ et $\betab:=(\beta_1,\beta_2)$ tels que $\beta_1=\alpha_1$ et $\beta_2^A=\alpha_2$. Comme $\betab$ est à coordonnées multiplicativement indépendantes, toute courbe $\Qr$-irréductible de $\Gm^2$ passant par $\betab$ n'est pas de torsion. D'après le théorème~\ref{Mthm} nous avons alors~:
$$
h(\betab)\geq \frac{1,\!2\cdot 10^{-16}}{\io_{\Qr}(\betab)}\max\left\{\log \io_{\Qr},\log 16\right\}^{-13}.
$$
On choisit maintenant :
$$
A:=\left[ \frac{2h(\alpha_2)}{h(\alpha_1)}\right]>\frac{2h(\alpha_2)}{h(\alpha_1)}-1\geq \frac{h(\alpha_2)}{h(\alpha_1)}
$$
en particulier
$$
h(\betab)\leq h(\beta_1)+h(\beta_2)= h(\alpha_1)+A^{-1}h(\alpha_2)\leq 2h(\alpha_1)\ ,
$$
et
$$
\io_{\Qr}(\betab)\leq 2[\Qr(\betab):\Qr]^{1/2}\leq 2(AD)^{1/2}\leq 2(2\frac{h(\alpha_2)}{h(\alpha_1)}D)^{1/2}.
$$
Nous obtenons alors
$$
2h(\alpha_1)\geq 
\frac{1,\!2\cdot 10^{-16}}{2(2\frac{h(\alpha_2)}{h(\alpha_1)}D)^{1/2}}
\max\left\{\log (2(AD)^{1/2}),\log 16\right\}^{-13}
$$
soit
\begin{equation}\label{eqCor}
\left(h(\alpha_1)h(\alpha_2)\right)^{1/2}\geq 
\frac{3.10^{-17}}{(2D)^{1/2}}\max\left\{\log (2(AD)^{1/2}),\log 16\right\}^{-13}.
\end{equation}
Minorons le membre de droite, nous avons :
$$
2(AD)^{1/2} \leq 2\left(\frac{2D}{h(\alpha_1)}\right)^{1/2}.
$$
Comme $\alphab$ est à coordonnées multiplicativement indépendantes,
$\alpha_1$ n'est pas une racine de l'unité et la version explicite du
théorème de Dobrowolski par Voutier~\cite{Vo} nous donne~:
$$
h(\alpha_1)\geq \frac{1}{4D} \left(\log (3D)\right)^{-3}.
$$
On en déduit
$$
\begin{disarray}{lcl}
\log \left(2(AD)^{1/2}\right)        & \leq &
 \log \left(4\sqrt{2}D(\log 3D)^{3/2}\right)\\[4mm]
                                        & \leq &
 2\log (3D).
\end{disarray}
$$
Comme $D\geq 2$, nous avons $2\log (3D)\geq \log 16$, et, en reportant ceci dans~(\ref{eqCor}) nous obtenons~:
$$
\left(h(\alpha_1)h(\alpha_2)\right)^{1/2}\geq 
\frac{3\cdot 2^{-13,\!5}\cdot 10^{-17}}{D^{1/2}}\left(\log (3D)\right)^{-13}.
$$
Ainsi $C(2)^{1/2}=2.10^{-21}$ et $\kappa(2)=13$.


\bigskip

\noindent Laboratoire de math\'ematiques Nicolas Oresme, CNRS UMR 6139\\
Universit\'e de Caen BP 5186\\
14032 Caen Cedex\\
FRANCE\\
E-mail~: pontreau@math.unicaen.fr\\


\begin{thebibliography}{WWW}
       \bibitem[Am-Da]{Am-Da} F. Amoroso et S. David, {\it Le probl\`eme de
Lehmer en dimension sup\'erieure}, J. Reine Angew. Math. 513 (1999), p. 145-179.
%
       \bibitem[Am-Da2]{Am-Da2} F. Amoroso et S. David, {\it Minoration de la hauteur normalisée d'une hypersurface}, Acta Arithmetica 92.4 (2000),
p. 339-366.
%
       \bibitem[Am-Da3]{Am-Da3} F. Amoroso et S. David, {\it Distribution des points de petite hauteur dans les groupes multiplicatifs}, Ann. Scuola Norm. Sup. Pisa Sci. Serie V Vol III Fasc. 2 (2004), p. 325-348.
%
       \bibitem[Am-Da4]{Am-Da4} F. Amoroso et S. David, {\it Minoration de la hauteur normalisée dans un tore}, Journal of the Inst. of Math. Jussieu (2003) {\bf 2}(3), p. 335-381.
%
       \bibitem[Am-Za]{Am-Za} F. Amoroso et U. Zannier, {\it A Relative Dobrowolski Lower Bound over Abelian Extensions}, Ann. Scuola Nom. Sup. Pisa Serie IV {\bf 29} (2000), p. 711-727.
%
        \bibitem[Bo-Va]{Bo-Va} E. Bombieri and J. Vaaler, {\it On
Siegel's Lemma}, Inv. math., 1983.
%
        \bibitem[Bi]{Bi} Y. Bilu, Math. Reviews MR 2000g:11058.
%
        \bibitem[Ch]{Ch} M. Chardin, {\it Une majoration de la fonction de Hilbert et ses conséquences pour l'interpolation algébrique.}, Bul. Soc. Math. France t. 117 (1988), p. 305-318.
%
          \bibitem[Da-Ph]{Da-Ph} S. David et P. Philippon, {\it Minorations des hauteurs normalisées des sous-variétés des tores}, Ann. Scuola Norm. Sup. Pisa Sci. (4) XXVIII (1999), pp. 489-543.
%
          \bibitem[Do]{Do} E. Dobrowolski, {\it On a question of Lehmer and the number of irreductible factors of a polynomial}, Acta Arithmetica 34 (1979),
p. 391-401.
%
	\bibitem[Le]{Le} D. H. Lehmer, {\it Factorization of certain cyclotomic functions}, Ann. Math. (2) 34, (1933), p. 461-479.
%
	\bibitem[P]{DEA} C. Pontreau, {\it Mémoire de DEA}, Université de Caen 2001
%
        \bibitem[Ro-Sc]{Ro-Sc} J. B. Rosser and L. Sh\oe nfeld, {\it Approximate formulas for some functions of prime numbers}, Ill. J. Math. t. {\bf 6}, pages 64-94, 1962.
%
           \bibitem[St-Va]{St-Va} T. Struppeck and J. Vaaler, {\it Inequalities for heights of algebrais subspaces and the Thue-Siegel principle}, in :
Analytic Number Theory, Boston, 1990, p. 494-527.
%
        \bibitem[Vo]{Vo} P. Voutier, {\it An effective lower bound for the height of algebraic numbers}, Acta Arithmetica 74, (1996), p. 81-95.
%
         \bibitem[Zh]{Zh} S. Zhang, {\it Positive line bundles on arithmetic
varieties}, J. Amer. Math. Soc. 8, (1995), p. 187-221.
%
\end{thebibliography}
\end{document}